\documentclass[preprint]{elsarticle}

\usepackage{amsmath,amssymb}
\usepackage{booktabs}
\usepackage{comment}
\usepackage{physics}
\usepackage[colorlinks,allcolors=blue]{hyperref}
\usepackage{cleveref} 
\usepackage{float}
\hypersetup{
    pdftitle={Fourier transform-based computation of two-electron four-center Coulomb repulsion integrals},
    pdfauthor={William Chen, Ury Segal, and Chen Greif}
}
\allowdisplaybreaks

\journal{Journal of Chemical Physics}

\begin{document}

\begin{frontmatter}

    \title{Fourier Transform-Based Computation of Two-Electron Four-Center Coulomb Repulsion Integrals}

    \author[1]{William Chen}
    \ead{bikachen@student.ubc.ca}
    \author[2]{Lior Silberman\fnref{fn2}}
    \ead{lior@math.ubc.ca}
    \author[3]{Ury Segal}
    \ead{ury.segal@gmail.com}
    \author[1]{Chen Greif\corref{cor1}\fnref{fn3}}
    \ead{greif@cs.ubc.ca}

    \cortext[cor1]{Corresponding author}
    \fntext[fn2]{The work of this author was funded by a Discovery Grant from the Natural Sciences and Engineering Research Council of Canada (NSERC) under grant number RGPIN-2019-03964.}
    \fntext[fn3]{The work of this author was funded by a Discovery Grant from the Natural Sciences and Engineering Research Council of Canada (NSERC) under grant number RGPIN-2023-05244.}
    \address[1]{The University of British Columbia, Department of Computer Science, Vancouver V6T 1Z4, Canada}
    \address[2]{The University of British Columbia, Department of Mathematics, Vancouver V6T 1Z2, Canada}
    \address[3]{QEndeavors Ltd.}

    \begin{abstract}
        We present a Fourier-based method for evaluating two-electron four-center Coulomb repulsion integrals over Slater-type orbitals. Spherical Fourier-Bessel coefficients of the convolution factors are precomputed and reused.  Since the spherical part of the integral can be evaluated directly, there remains only a one-dimensional quadrature problem in the radial direction.  For the ten reference integrals reported here, the largest absolute discrepancy  is $1.7\times10^{-9}$. A preliminary resident-cache benchmark computes approximately $4.2\times10^5$ unique integrals per second on an NVIDIA V100 at the contraction stage.
    \end{abstract}

    \begin{keyword}
        four-center integrals \sep Slater-type orbitals \sep Fourier transforms
        \sep spherical harmonics \sep analytic angular integration \sep GPU
    \end{keyword}

\end{frontmatter}

\section{Introduction}
\label{sec:introduction}

Fast ab initio calculations of molecular properties require efficient
evaluation of electron repulsion integrals (ERIs). In orbital-based
methods, molecular orbitals are expressed as linear combinations of
basis functions
\cite{szaboModernQuantumChemistry1996,cookHandbookComputationalQuantum2005}.
Each ERI describes the Coulomb interaction between two orbital-pair
densities and involves four basis functions. Slater-type orbitals
(STOs)~\cite{slaterSelfConsistentField1928a} are a possible choice. They
give rise to accurate computations using fewer basis functions
compared to competing approaches, such as Gaussian-type orbitals, but
the resulting multicenter integrals are known to be difficult to evaluate.

The ERI can be viewed as a convolution of two functions, each of which is the product of two STOs, mediated by the Coulomb potential as a kernel.  It is therefore natural to replace the six-dimensional integral with the three-dimensional integration of the Fourier transforms of the two factors against the Fourier transform of the potential.  A primary
computational challenge is to represent these transforms compactly and
reuse them across the many integrals sharing the same orbital pairs.

Partial-wave analysis has long been central to this approach. Trivedi
and Steinborn derived a two-center transform of reduced Bessel functions
explicitly suited to partial-wave decomposition and proposed reusing
previously computed transforms \cite{trivedi1983fourier}. Grotendorst
and Steinborn subsequently developed a general multicenter treatment
using  $B$ functions, with a three-dimensional numerical
integration remaining for two-electron integrals
\cite{grotendorst1988numerical}. Direct STO formulations include the
shift-operator construction of Niehaus, L\'opez, and Rico
\cite{niehaus2008efficient} and the one-dimensional, non-oscillatory
hypergeometric representation of Vukovi\'c and Dmitrovi\'c
\cite{vukovic2010extremely}. These developments provide analytic
representations of pair transforms and multicenter integrals.

Our focus
is on the organization of reusable local pair data and angular contractions
for the evaluation of large batches of ERIs on GPUs.
We developed a localized, reusable partial-wave formulation. A
pair-specific Becke partition suppresses the off-center nuclear
structure, allowing each component to be expanded about its own center.
The center phases are kept outside these expansions and analytically
integrated during contraction. All local harmonic channels are stored
at each radial momentum node, and the resulting finite angular
couplings are evaluated in GPU batches.

The remainder of the paper is organized as follows.
Section~\ref{sec:fourier_formulation} presents the Fourier formulation
of the ERIs. In Section~\ref{sec:local_representation} we construct the
localized pair representation and its reusable coefficients.
In Section~\ref{sec:angular_contraction} the analytic angular
contraction is derived. In Section~\ref{sec:implementation} we describe the cache layout
and computational organization. Numerical results are presented in
Section~\ref{sec:results}, followed by brief concluding remarks in Section~\ref{sec:conclusions}.

\section{Fourier formulation}
\label{sec:fourier_formulation}

Let $\{\phi_\mu\}_{\mu=1}^{N}$ be real, normalized STOs,
\begin{equation*}
    \phi_\mu(\mathbf r)
    =\mathcal N_\mu r_\mu^{\,n_\mu-1}e^{-\zeta_\mu r_\mu}
    Y^{\mathrm R}_{\ell_\mu m_\mu}\qty(\frac{{\mathbf r}_\mu}{r_\mu}),
\end{equation*}
where  $\mathbf r_\mu=\mathbf r-\mathbf A_\mu,$
    $r_\mu=|\mathbf r_\mu|$, $n_\mu$ is the principal index, $\zeta_\mu>0$, and
$Y^{\mathrm R}_{\ell m}$ denotes an orthonormal basis of real spherical harmonics.
For $p=(\mu,\nu)$ with $\mu\leq\nu$, set
$\rho_p=\phi_\mu\phi_\nu$. There are $M=N(N+1)/2$ such pairs. The ERI
associated with $p=(\mu,\nu)$ and $q=(\lambda,\sigma)$ is
\begin{equation*}
    C_{pq}=(\mu\nu\mid\lambda\sigma)
    =\iint_{\mathbb R^3\times\mathbb R^3}
    \frac{\rho_p(\mathbf r)\rho_q(\mathbf r')}
    {|\mathbf r-\mathbf r'|}\,d\mathbf r\,d\mathbf r'.
\end{equation*}

Denoting the Fourier transform by
\begin{equation*}
    F_p(\mathbf k)=\int_{\mathbb R^3}\rho_p(\mathbf r)
    e^{-i\mathbf k\cdot\mathbf r}\,d\mathbf r,
\end{equation*}
the convolution theorem and Parseval's identity give
\begin{equation*}
    C_{pq}=\frac{1}{(2\pi)^3}\int_{\mathbb R^3}
    \frac{4\pi}{|\mathbf k|^2}
    \overline{F_p(\mathbf k)}F_q(\mathbf k)\,d\mathbf k.
\end{equation*}
Writing $\mathbf k=k\boldsymbol\omega$, with
$\boldsymbol\omega\in S^2$ and $k\geq0$, we have the volume element
$d\mathbf k=k^2\,dk\,d\boldsymbol\omega$. The spherical Jacobian
cancels the $k^{-2}$ Coulomb multiplier, giving
\begin{equation}
    C_{pq}=\frac{1}{2\pi^2}\int_0^\infty\int_{S^2}
    \overline{F_p(k\boldsymbol\omega)}F_q(k\boldsymbol\omega)
    \,d\boldsymbol\omega\,dk.
    \label{eq:fourier_eri}
\end{equation}

A common-origin expansion
$F_p(k\boldsymbol\omega)=\sum_{\ell m}
    f_{p,\ell m}(k)Y_{\ell m}(\boldsymbol\omega)$ consequently gives
\begin{equation*}
    C_{pq}=\frac{1}{2\pi^2}\int_0^\infty
    \sum_{\ell m}\overline{f_{p,\ell m}(k)}f_{q,\ell m}(k)\,dk.
\end{equation*}
Unfortunately, using this will result in a larger cache size,
described in Section~\ref{sec:implementation}, because a translated spherical density has the transform
$\widehat g(k)e^{-ik\boldsymbol\omega\cdot\mathbf X}$, whose Rayleigh
coefficients contain $j_\ell(k|\mathbf X|)$ through degrees of order
$k|\mathbf X|$. For a molecule of radius $R_{\rm sys}$ relative to the
common origin, the global cutoff can therefore grow on the scale of
$kR_{\rm sys}$. We instead retain translations explicitly and expand
only locally-centered amplitudes.

\section{Localized representation}
\label{sec:local_representation}

Consider a pair with distinct centers $\mathbf A$ and $\mathbf B$.
An expansion about $\mathbf A$ encounters the displaced nuclear cusp on
shells passing through $\mathbf B$, where a radial feature becomes a
localized angular feature. We use a two-center Becke partition
\cite{becke1988multicenter} to suppress that feature in the
$\mathbf A$-centered component:
\begin{equation*}
    \rho_p=\rho_p^A+\rho_p^B,
    \qquad
    \rho_p^A=w_A\rho_p,
    \qquad
    \rho_p^B=w_B\rho_p.
\end{equation*}
Let $d_{AB}=|\mathbf A-\mathbf B|$ and define
\begin{equation*}
    \begin{aligned}
        t(\mathbf r)
                       & =\frac{|\mathbf r-\mathbf A|-|\mathbf r-\mathbf B|}{d_{AB}},
                       & s(t)                                                         & =\frac32t-\frac12t^3, \\
        w_A(\mathbf r) & =\frac12\bigl[1-s(s(s(t(\mathbf r))))\bigr],
                       & w_B(\mathbf r)                                               & =1-w_A(\mathbf r).
    \end{aligned}
\end{equation*}
The partition is exact before truncation. More importantly, endpoint
flatness gives $w_A=\mathcal O((r_B/d_{AB})^8)$ near $\mathbf B$, where
$r_B=|\mathbf r-\mathbf B|$, and conversely for $w_B$ near $\mathbf A$.
Each component therefore suppresses the other center's cusp to high
order while retaining its own in local radial coordinates. Same-center
pairs require no split.

Let $\mathcal C(p)$ denote the resulting one or two component centers.
For $X\in\mathcal C(p)$, write
$\mathbf r=\mathbf X+r\boldsymbol\Omega$ and project the component onto
real spherical harmonics:
\begin{equation*}
    \rho^X_{p,\ell m}(r)
    =\int_{S^2}\rho_p^X(\mathbf X+r\boldsymbol\Omega)
    Y^{\mathrm R}_{\ell m}(\boldsymbol\Omega)
    \,d\boldsymbol\Omega.
\end{equation*}
The reusable radial coefficients are
\begin{equation}
    R^X_{p,\ell m}(k)
    =\int_0^\infty r^2j_\ell(kr)\rho^X_{p,\ell m}(r)\,dr,
    \label{eq:radial_coefficients}
\end{equation}
where $j_\ell$ is the spherical Bessel function of the first kind.
These coefficients are real and are evaluated numerically for
$0\leq\ell\leq L$, $\ -\ell\leq m\leq\ell$, where $L=L_{\max}$.

The Fourier transform of a degree-$\ell$ spherical-harmonic component
has the same angular harmonic, multiplied by $4\pi(-i)^\ell$ and its
radial spherical-Bessel transform. Applying this identity to each
component and retaining degrees through $L$ defines
\begin{equation}
    \begin{aligned}
        F_p^{(L)}(k\boldsymbol\omega)
         & =4\pi\sum_{X\in\mathcal C(p)}
        e^{-ik\boldsymbol\omega\cdot\mathbf X}
        G_p^X(k,\boldsymbol\omega),                \\
        G_p^X(k,\boldsymbol\omega)
         & =\sum_{\ell=0}^{L}\sum_{m=-\ell}^{\ell}
        (-i)^\ell R^X_{p,\ell m}(k)
        Y^{\mathrm R}_{\ell m}(\boldsymbol\omega).
    \end{aligned}
        \label{eq:local_fourier}
\end{equation}
Only $G_p^X$ is represented by a finite angular expansion: the
translation phase is neither expanded nor sampled. Localization
reduces the degree needed for off-center structure.

\section{Analytic angular contraction}
\label{sec:angular_contraction}

Our goal is to evaluate the angular integral in
\eqref{eq:fourier_eri} directly from the cached coefficients, leaving
a one-dimensional momentum quadrature. Substituting the local
representation \eqref{eq:local_fourier} reduces this task to
component-center pairs. For each pair $(X,Y)$, write
$\mathbf D_{XY}=\mathbf X-\mathbf Y$ and $D_{XY}=|\mathbf D_{XY}|$.
The required contribution is
\begin{equation*}
    \mathcal A_{pX,qY}(k)
    =
    \int_{S^2}
    e^{ik\boldsymbol\omega\cdot\mathbf D_{XY}}
    \overline{G_p^X(k,\boldsymbol\omega)}
    G_q^Y(k,\boldsymbol\omega)
    \,d\boldsymbol\omega.
\end{equation*}
Summing up these contributions and multiplying by $(4\pi)^2$ gives
the angular integral of
$\overline{F_p^{(L)}}F_q^{(L)}$.

For the derivation, we express the local amplitudes in complex
spherical harmonics:
\begin{equation*}
    G_p^X(k,\boldsymbol\omega)
    =
    \sum_{\ell=0}^{L}\sum_{m=-\ell}^{\ell}
    u^X_{p,\ell m}(k)Y_{\ell m}(\boldsymbol\omega),
    \qquad
    u^X_{p,\ell m}(k)=(-i)^\ell c^X_{p,\ell m}(k),
\end{equation*}
where $c^X_{p,\ell m}$ is the unitary conversion of the cached real
coefficients $R^X_{p,\ell m}$.

For distinct centers, we rotate the coordinates to align
$\mathbf D_{XY}$ with the positive $z$-axis. Let $Q$ denote this
rotation and $\boldsymbol\eta=Q\boldsymbol\omega$ the rotated
direction. Both amplitudes are expressed in this same frame using
the Wigner rotation matrices:
\begin{equation*}
    \begin{aligned}
        Y_{\ell s}(Q^T\boldsymbol\eta)
         & =
        \sum_{m=-\ell}^{\ell}
        S^{(\ell)}_{ms}(Q)Y_{\ell m}(\boldsymbol\eta), \\
        \widetilde u^X_{p,\ell m}(k)
         & =
        \sum_{s=-\ell}^{\ell}
        S^{(\ell)}_{ms}(Q)u^X_{p,\ell s}(k).
    \end{aligned}
\end{equation*}

In this frame, the phase becomes $e^{izt}$, with $z=kD_{XY}$ and
$t$ the cosine of the polar angle. Azimuthal integration retains
equal magnetic orders, giving the translation blocks
\begin{equation*}
    T^m_{\ell\ell'}(z)
    =
    2\pi N_{\ell a}N_{\ell'a}
    \int_{-1}^{1}
    P_\ell^a(t)P_{\ell'}^a(t)e^{izt}\,dt,
    \qquad a=|m|,
\end{equation*}
where
$N_{\ell a}=\sqrt{(2\ell+1)(\ell-a)!/[4\pi(\ell+a)!]}$
is the spherical-harmonic normalization.
For each $m$, the degrees $\ell,\ell'=|m|,\ldots,L$ form a block
of order $L-|m|+1$.

To evaluate these blocks, we precompute the normalized
Legendre-product expansions
\begin{equation*}
    N_{\ell a}N_{\ell'a}P_\ell^a(t)P_{\ell'}^a(t)
    =
    \sum_{n=0}^{\ell+\ell'}
    b^a_{n,\ell\ell'}P_n(t).
\end{equation*}
Combining this expansion with
$\int_{-1}^{1}P_n(t)e^{izt}\,dt=2i^n j_n(z)$ gives
\begin{equation*}
    T^m_{\ell\ell'}(z)
    =
    4\pi\sum_{n=0}^{\ell+\ell'}
    b^{|m|}_{n,\ell\ell'}i^n j_n(z).
\end{equation*}
The coefficients depend only on angular indices. The selection
rules restrict their support to
$|\ell-\ell'|\leq n\leq\ell+\ell'$ with $\ell+\ell'+n$ even.
Thus Bessel orders through $2L$ suffice for every value of
$kD_{XY}$.

Collect the rotated coefficients of each magnetic $(\widetilde u^X_{p,\ell m}(k))_{\ell=|m|}^{L}$ order into
$\widetilde{\mathbf u}_{p,m}^{X}(k)$.
The component contribution is then
\begin{equation*}
    \mathcal A_{pX,qY}(k)
    =
    \sum_{m=-L}^{L}
    \widetilde{\mathbf u}_{p,m}^{X\dagger}(k)
    T^m(kD_{XY})
    \widetilde{\mathbf u}_{q,m}^{Y}(k).
\end{equation*}
For coincident centers, the original frame gives $T^m(0)=I$, the identity, 
and the contraction reduces to an ordinary coefficient
inner product.

Finally, using quadrature nodes and weights
$\{k_i,v_i\}_{i=1}^{N_k}$ on $[0,K_{\max}]$, we obtain
\begin{equation*}
    \widehat C_{pq}
    =
    8\sum_{i=1}^{N_k}v_i
    \sum_{X\in\mathcal C(p)}
    \sum_{Y\in\mathcal C(q)}
    \operatorname{Re}\mathcal A_{pX,qY}(k_i).
\end{equation*}
The factor $8=(4\pi)^2/(2\pi^2)$ combines the Fourier and Coulomb
prefactors. The angular integration is exact for the retained
harmonic representation.

\section{Implementation and storage}
\label{sec:implementation}
To generate the coefficients in \eqref{eq:radial_coefficients}, the
angular projections are evaluated using 48 Gauss--Legendre nodes in
$u=\cos\theta$ and 96 uniformly spaced azimuthal nodes.
We note that Lebedev quadrature is more efficient, as we observed in \cite{jcc.70472}, but we choose to use Gauss-Legendre for simplicity in terms of implementation; this part of the computation entails a relatively insignificant computational effort anyway.

The radial integrals use 24 panels with 16 Gauss--Legendre nodes per
panel. The grid is shared within each center pairs for all orbitals:
\[
R^X_{p,\ell m}(k_i)
\approx
\sum_{a=1}^{24 \times 16}
\underbrace{w_a^{(r)}r_a^2j_\ell(k_i r_a)}_{\text{reused across orbital pairs}}
\;\rho^X_{p,\ell m}(r_a).
\]
We also enforced a panel boundary at $r=d_{AB}$ to ensure each side is approximated separately, as discussed in Section \ref{sec:local_representation}.
We determine the  choices of numerical parameters empirically and note that they can be easily adjusted to fit the application needs.
For each center pair, we choose the radial cutoff to be
$r_{\max}=d_{AB}+\tau/\beta_{\min}$, with $\beta_{\min}$ being the smallest sum
$\zeta_\mu+\zeta_\nu$ among its orbital pairs.
For $p=(\mu,\nu)$, the STO decay and $|j_\ell(kr)|\le1$ give the
following bound on the tail:
\begin{equation*}
    \left|
    \int_{r_{\max}}^\infty
    r^2 j_\ell(kr)\rho^X_{p,\ell m}(r)\,dr
    \right|
    \leq C_p(1+\tau)^{n_\mu+n_\nu}e^{-\tau},
    \label{eq:radial_tail_bound}
\end{equation*}
where the constant \(C_p\) depends on the orbital parameters
and geometry, but is independent of $k$ and $\tau$. 
We choose \(\tau=30\) for the computation demonstrated in Section \ref{sec:results}.

The persistent cache stores
$R^X_{p,\ell m}(k_i)$ as real double-precision values in blocks indexed
by component, together with pair and center maps. Every block retains
all $H=(L+1)^2$ channels at every momentum node; complex coefficients
and rotations are formed for contraction. With
$N_{\rm comp}=\sum_p|\mathcal C(p)|$ and
$M\leq N_{\rm comp}\leq2M$, the coefficient storage in bytes is
\begin{equation}
    S_{\rm pair}=8N_{\rm comp}N_k(L+1)^2\ 
    =\mathcal O\!\left(N^2N_k(L+1)^2\right).
    \label{eq:storage}
\end{equation}
Thus the cache is quadratic in the number of basis functions at fixed
$N_k$ and $L$. Keeping translations outside the cache avoids a factor
set by $K_{\max}R_{\rm sys}$ in the stored angular degree. The dense
layout also avoids irregular per-pair channel lists. Device residency
is nonetheless limited by this quadratic cache and the additional
contraction workspaces.

Contractions are grouped by component-center pair so that rotations and
translation matrices are reused, with spherical-Bessel values shared
at repeated distances. The inner exact angular contraction exposes batched matrix operations
rather than repeated angular quadratures. For a single component-pair
contribution at one momentum node, the direct bilinear form involves
\begin{equation*}
    \sum_{m=-L}^{L}(L-|m|+1)^2
    =\frac{(L+1)(2L^2+4L+3)}{3}
\end{equation*}
matrix entries. Thus, once rotations and coupling matrices are available,
this gives a direct evaluation cost of
$\mathcal O(N_k(L+1)^3)$ per component-pair contribution, with further
reuse available across a batch. Each ERI accumulates at most four
component contributions.

Pair preprocessing is amortized over the integrals that share each
pair. Explicitly forming all unique ERIs produces
$M(M+1)/2$ output values, which scales as $\mathcal O(N^4)$.

\section{Results}
\label{sec:results}

Table~\ref{tab:results} reports the four-center test cases
from~\cite{jcc.70472}. The reference values are obtained from the
FORTRAN code described in~\cite{berluMulticentreTwoelectronCoulomb2004}.
The four orbitals in each entry are associated, in order, with centers
$\mathbf c_1$, $\mathbf c_2$, $\mathbf c_3$, and $\mathbf c_4$ given in
the caption. The comparison uses the orbital labeling of the cited
test cases.

The examples include $s$-orbital products, mixed $s$--$p$ products, and
several $d$-orbital combinations. The reported absolute discrepancies
are at or below $10^{-10}$ for seven entries and between
$1.5\times10^{-9}$ and $1.7\times10^{-9}$ for the remaining three.
The largest reported discrepancy is $1.7\times10^{-9}$, occurring in
two of the $d$-orbital examples. The signed mixed-orbital values and the
different magnetic orders provide further evidence of the validity of the angular factors in
addition to the radial dependence.

\begin{table}[htbp]
    \centering
    \small
    \setlength{\tabcolsep}{5pt}
    \begin{tabular}{|c|r|r|r|} \hline
        \textbf{Integral $(\Phi_1\Phi_2\mid\Phi_3\Phi_4)$} & \textbf{Result}  & \textbf{Ref Result} & \textbf{Abs. Error}     \\ \hline
        $(1s\,1s\mid1s\,1s)$                               & $0.3111837840  $ & $ 0.3111837840   $  & $< 1.0 \times 10^{-10}$ \\ \hline
        $(1s\,1s\mid2s\,1s)$                               & $0.2622985772  $ & $ 0.2622985772   $  & $< 1.0 \times 10^{-10}$ \\ \hline
        $(2s\,2s\mid2s\,2s)$                               & $0.2849043959  $ & $ 0.2849043960   $  & $1.0\times10^{-10}$     \\ \hline
        $(1s\,1s\mid2s\,2p_{-1})$                          & $-0.1226144752 $ & $ -0.1226144752  $  & $< 1.0 \times 10^{-10}$ \\ \hline
        $(1s\,2s\mid2s\,2p_{-1})$                          & $-0.1006677390 $ & $ -0.1006677390  $  & $< 1.0 \times 10^{-10}$ \\ \hline
        $(2s\,2p_{-1}\mid1s\,1s)$                          & $-0.0038018631 $ & $ -0.0038018631  $  & $< 1.0 \times 10^{-10}$ \\ \hline
        $(3d_{-1}\,3d_{-1}\mid3d_{-1}\,3d_{-1})$           & $0.1740122879  $ & $ 0.1740122896   $  & $1.7\times10^{-9}~$     \\ \hline
        $(3d_{+2}\,3d_{+1}\mid3d_{+2}\,3d_{+1})$           & $0.0026966172  $ & $ 0.0026966172   $  & $< 1.0 \times 10^{-10}$ \\ \hline
        $(3d_{+1}\,3d_{+1}\mid3d_{+2}\,3d_{+2})$           & $0.1531260148  $ & $ 0.1531260163   $  & $1.5\times10^{-9}~$     \\ \hline
        $(3d_{+1}\,3d_{+1}\mid3d_{+1}\,3d_{+1})$           & $0.1740122879  $ & $ 0.1740122896   $  & $1.7\times10^{-9}~$     \\
        \hline
    \end{tabular}
    \caption{Four-center STO electron-repulsion integrals for centers
        $\mathbf c_1=(0,0,0)$, $\mathbf c_2=(1,0,0)$,
        $\mathbf c_3=(2,0,0)$, and $\mathbf c_4=(3,0,0)$ 
        with $L_{\max}=22$, $K_{\max}=40$, and a momentum quadrature
        of 24 panels with 16 Gauss--Legendre nodes per panel.
    }
    \label{tab:results}
\end{table}

\begin{figure}[H]
    \centering
     \includegraphics[scale=0.4]{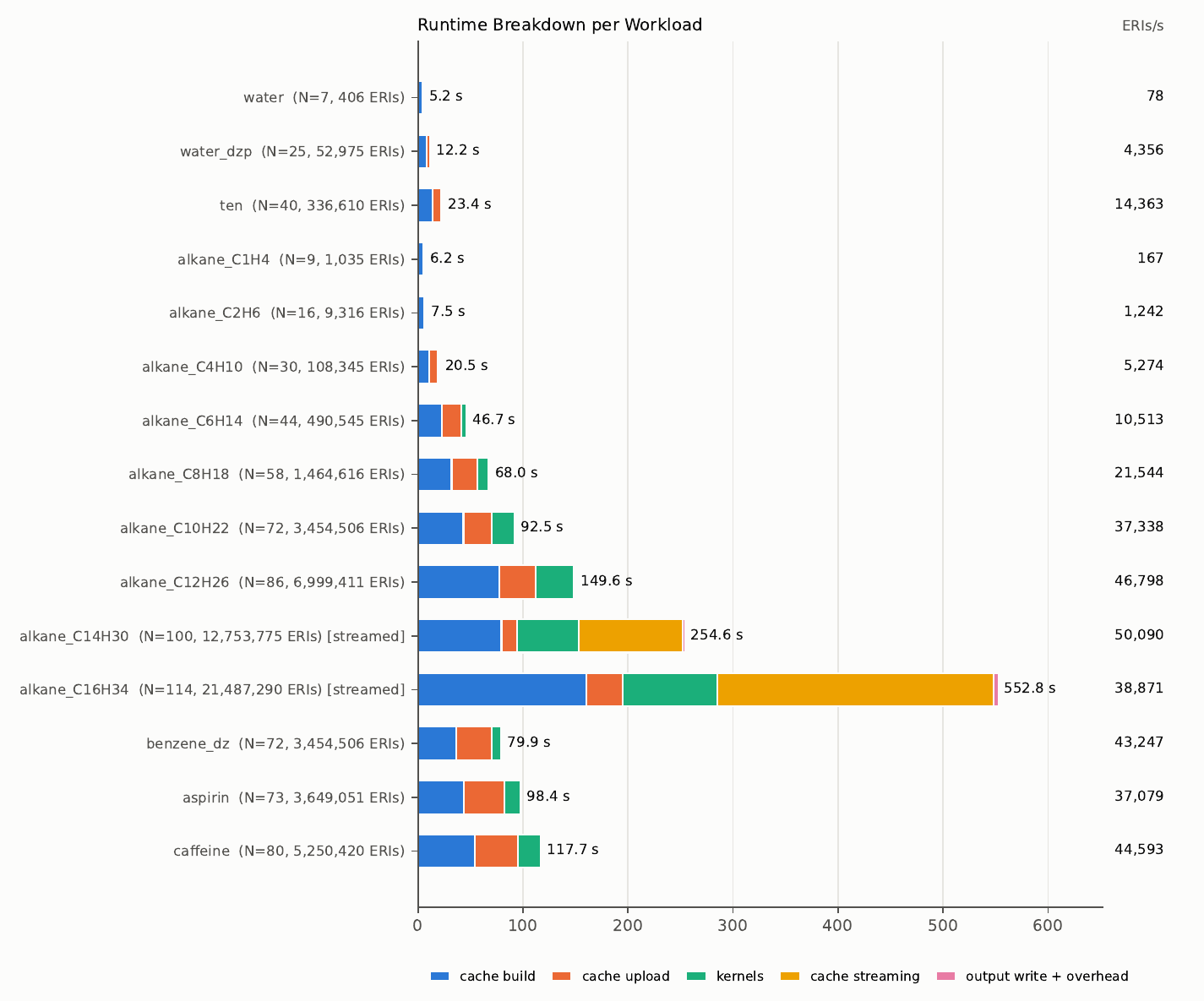}
     \label{fig:runtime}
     \caption{Total runtime breakdown per workload using the parameters of \Cref{tab:results}.}
\end{figure}

Figure \ref{fig:runtime} shows the end-to-end runtime for various tests, which is the wall clock time
 to compute all ERIs of a molecule. This includes pair cache generation, IO time, memory transfer, and additional computational overhead. The time taken for
the cache build depends on the geometry of the molecule and the number of basis functions. If the
total memory of the workspace used by the contraction stage exceeds the GPU memory,
then the workload must be split and streamed in smaller parts, resulting in more time spent on memory transfer. Overall, the diagram shows that we sustain about 40,000 ERIs per second end-to-end.

A contraction benchmark reaches
approximately $4.2\times10^5$ unique ERIs per second on an NVIDIA V100.
This is a contraction-stage rate, not an end-to-end timing including
cache generation. The throughput measures the reuse of precomputed pair
coefficients during ERI assembly. Together with the storage formula
\eqref{eq:storage}, it characterizes the contraction rate and coefficient
storage of the implementation.

\section{Conclusions}
\label{sec:conclusions}

The localized Fourier formulation for evaluating
four-center Coulomb repulsion integrals over STOs that we have presented separates
pair preprocessing from ERI contraction. A two-center Becke partition
produces locally centered components, and their radial spherical-Bessel
coefficients provide reusable data for all integrals involving the
corresponding orbital pairs.

We take advantage of the attractive computational features of our GPU environment, and the results are promising. The small absolute discrepancies and the large number of unique ERIs that are computed per second suggest that further
study of the accuracy--storage--throughput balance through systematic
refinement of the local representation and momentum quadrature may be beneficial.

\bibliographystyle{elsarticle-num}
\bibliography{references}

\end{document}